\documentclass[amstex,12pt,russian,amssymb]{article}

\usepackage{mathtext}
\usepackage[cp1251]{inputenc}
\usepackage[T2A]{fontenc}
\usepackage[russian]{babel}
\usepackage[dvips]{graphicx}
\usepackage{amsmath}
\usepackage{amssymb}
\usepackage{amsxtra}
\usepackage{latexsym}
\usepackage{ifthen}

\begin{document}

\newcounter{lemma}
\newcommand{\lemma}{\par \refstepcounter{lemma}%
{\bf Лемма \arabic{lemma}.}}

\newcounter{corollary}
\newcommand{\corollary}{\par \refstepcounter{corollary}%
{\bf Следствие \arabic{corollary}.}}

\newcounter{remark}
\newcommand{\remark}{\par \refstepcounter{remark}%
{\bf Замечание \arabic{remark}.}}

\newcounter{theorem}
\newcommand{\theorem}{\par \refstepcounter{theorem}%
{\bf Теорема \arabic{theorem}.}}

\newcounter{proposition}
\newcommand{\proposition}{\par \refstepcounter{proposition}%
{\bf Предложение \arabic{proposition}.}}

\renewcommand{\refname}{\centerline{\bf Список литературы}}

\newcommand{\proof}{{\it Доказательство.\,\,}}

{\bf Е.А. Севостьянов} (Институт прикладной математики и механики
НАН Украины, Донецк, Украина)

{\bf E.A. Sevost'yanov} (Institute of Applied Mathematics and
Mechanics of NAS of Ukraine, Donetsk, Ukraine)

{\bf О равностепенной непрерывности открытых дискретных отображений
классов Орлича--Соболева}

{\bf Equicontinuity of the family of the open discrete
Orlicz--Sobolev mappings}

\medskip
{\small The paper is devoted to the study of mappings with
non--bounded characteristics of qua\-si\-con\-for\-ma\-li\-ty. We
investigate the interconnection between the classes of the
so--called ring $Q$-mappings and lower $Q$-mappings. It is proved
that open discrete lower ring $Q$-mappings are ring
$Q^{n-1}$-mappings at fixed point. As consequence we obtain the
equicontinuity of the class of the open discrete Orlicz--Sobolev
mappings with finite distortion at $n\ge 3.$}

\section{Введение}

\medskip{} Данная работа посвящена изучению свойства равностепенной
непрерывности одного подвида отображений с конечным искажением,
активно изучаемых последнее время (см. \cite{IM}). Как было показано
в одной из совместных работ автора, семейства гомеоморфизмов класса
Орлича--Соболева с конечным искажением являются нормальными
(равностепенно непрерывными) при определённых дополнительных
условиях на характеристику квазиконформности отображений и
количество выпускаемых этими отображениями значений (см., напр.,
\cite[теоремы~7 и 9, следствие~13]{KRSS}). Однако, как оказалось,
требование гомеоморфности в формулировке упомянутых результатов в
известном смысле не является принципиальным, поскольку, как будет
показано в настоящей работе, при некоторых (довольно естественных)
условиях на семейство отображений, условие гомеоморфности можно
отбросить и заменить его требованием, что каждое отображение
является лишь открытым и дискретным. При этом, дополнительно
налагается условие ограниченности функции кратности рассматриваемого
семейства отображений. В данной работе будут сформулированы и
подробно доказаны результаты подобного характера.

\medskip
Как и в случае гомеоморфизмов, исследование открытых дискретных
отображений классов Орлича--Соболева опирается на их связь с так
называемыми нижними и кольцевыми $Q$-отоб\-ра\-же\-ни\-я\-ми, в
связи с чем в работе развивается параллельная вспомогательная теория
их исследования. Ниже мы приведём наиболее часто использующиеся в
статье понятия, а также дадим формулировки основных результатов.

\medskip{} Всюду далее $D$ -- область в ${\Bbb R}^n,$ $n\ge 2,$ $m$ -- мера
Лебега в ${\Bbb R}^n$ и ${\rm dist\,}(A,B)$ -- евклидово расстояние
между множествами $A$ и $B$ в ${\Bbb R}^n.$ Запись $f:D\rightarrow
{\Bbb R}^n$ предполагает, что отображение $f$ непрерывно в $D.$ В
дальнейшем ${\mathcal H}^k$ -- нормированная $k$-мерная мера
Хаусдорфа в ${\Bbb R}^n,$ $1\le k\le n,$
$$B(x_0, r)=\left\{x\in{\Bbb R}^n: |x-x_0|< r\right\}\,,\quad {\Bbb B}^n
:= B(0, 1)\,,$$ $$S(x_0,r) = \{ x\,\in\,{\Bbb R}^n :
|x-x_0|=r\}\,,\quad{\Bbb S}^{n-1}:=S(0, 1)\,,$$
$$A(r_1,r_2,x_0)=\{ x\,\in\,{\Bbb R}^n : r_1<|x-x_0|<r_2\}\,,$$
$\omega_{n-1}$ обозначает площадь единичной сферы ${\Bbb S}^{n-1}$ в
${\Bbb R}^n,$ $\Omega_{n}$ -- объём единичного шара ${\Bbb B}^{n}$ в
${\Bbb R}^n.$ Напомним, что отображение $f:D\rightarrow {\Bbb R}^n$
называется {\it открытым,} если множество $f(E)$ открыто в ${\Bbb
R}^n$ для каждого открытого $E\subset D.$ Отображение
$f:D\rightarrow {\Bbb R}^n$ называется {\it дискретным,} если полный
прообраз $f^{-1}\left(y\right)$ каждой точки $y\in{\Bbb R}^n$
состоит только из изолированных точек. В дальнейшем $J(x, f)={\rm
det}\, f^{\,\prime}(x)$ -- {\it якобиан отображения} $f$ в точке
$x,$ где $f^{\,\prime}(x)$ -- {\it матрица Якоби} отображения $f$ в
точке $x.$

\medskip{} Пусть $\varphi:[0,\infty)\rightarrow[0,\infty)$ --
неубывающая функция, $f$ -- локально интегрируемая вектор-функция
$n$ вещественных переменных $x_1,\ldots,x_n,$ $f=(f_1,\ldots,f_m),$
$f_i\in W_{loc}^{1,1},$ $i=1,\ldots,m.$ Будем говорить, что
$f:D\rightarrow {\Bbb R}^n$ принадлежит классу
$W^{1,\varphi}_{loc},$ пишем $f\in W^{1,\varphi}_{loc},$ если
%
$$\int\limits_{G}\varphi\left(|\nabla
f(x)|\right)\,dm(x)<\infty\,,$$
%
для любой компактной подобласти $G\subset D,$ где $|\nabla
f(x)|=\sqrt{\sum\limits_{i=1}^m\sum\limits_{j=1}^n\left(\frac{\partial
f_i}{\partial x_j}\right)^2}.$ Класс $W^{1,\varphi}_{loc}$
называется классом {\it Орлича--Соболева}.

\medskip
Отображение $f:D\rightarrow {\Bbb R}^n$ называется {\it отображением
с конечным искажением}, пишем $f\in FD,$ если $f\in
W_{loc}^{1,1}(D)$ и для некоторой функции $K(x): D\rightarrow
[1,\infty)$ выполнено условие
$\Vert f^{\,\prime}\left(x\right) \Vert^{n}\le K(x)\cdot |J(x,f)|$
при почти всех $x\in D$ (см. \cite[п.~6.3, гл.~VI]{IM}. Для
отображений с конечным искажением корректно определена и почти всюду
конечна так называемая {\it внешняя дилатация} $K_O(x,f)$
отображения $f$ в точке $x,$ определяемая соотношением
\begin{equation}\label{eq0.1.1A}
K_O(x,f)\quad =\quad\left\{
\begin{array}{rr}
\frac{\Vert f^{\,\prime}(x)\Vert^n}{|J(x, f)|}, & J(x,f)\ne 0,\\
1,  &  f^{\,\prime}(x)=0, \\
\infty, & \text{в\,\,остальных\,\,случаях}
\end{array}
\right.\,.
\end{equation}
Для отображения $f:D\,\rightarrow\,{\Bbb R}^n,$ множества $E\subset
D$ и $y\,\in\,{\Bbb R}^n,$  определим {\it функцию кратности $N(y,
f, E)$} как число прообразов точки $y$ во множестве $E,$ т.е.
\begin{equation}\label{eq1.7A}
N(y, f, E)\,=\,{\rm card}\,\left\{x\in E: f(x)=y\right\}\,,
%
N(f, E)\,=\,\sup\limits_{y\in{\Bbb R}^n}\,N(y, f, E)\,.
\end{equation}

\medskip{}\label{p1.4} Пусть $G$ -- открытое множество в ${\Bbb R}^n$ и $I=\{x\in{\Bbb
R}^n:a_i<x_i<b_i,i=1,\ldots,n\}$ -- открытый $n$-мерный интервал.
Говорят, что отображение $f:I\rightarrow{\Bbb R}^n$ {\it принадлежит
классу $ACL$} ({\it абсолютно непрерывно на линиях}), если $f$
абсолютно непрерывно на почти всех линейных сегментах в $I,$
параллельных координатным осям. Говорят, что отображение
$f:G\rightarrow{\Bbb R}^n$ {\it принадлежит классу $ACL$} в $G,$
если сужение $f|_I$ принадлежит классу $ACL$ для каждого интервала
$I,$ $\overline{I}\subset G.$ Пусть $U$ -- открытое множество в
${\Bbb R}^n.$ Будем говорить, что $f\in ACL^p(U),$ если $f\in
ACL(U)$ и все частные производные $\frac{\partial f_i}{\partial
x_j}(x)$ по каждой из переменных $x_j$ в точке $x,$ $1\le i\le n,$
$1\le j\le n,$ локально интегрируемы в $U$ в степени $p.$ Пусть
$p\ge 1,$ $U$ -- открытое множество в ${\Bbb R}^n,$ тогда
$W_{loc}^{1, p}(U)=ACL^p(U),$
см. \cite[теоремы~1 и 2, п.~1.1.3, $\S\,$1.1, гл.~I]{Ma}.

\medskip{}
{\it Конденсатором} называют пару $E=\left(A,\,C\right),$ где $A$ --
открытое множество в ${\Bbb R}^n,$ а $C$ -- компактное подмножество
$A.$ {\it Ёмкостью} конденсатора $E$ называется следующая величина:
\begin{equation}\label{eq1.1AB} {\rm cap}\,E={\rm
cap}\,\left(A,\,C\right)= \inf\limits_{u\in W_0(E)}\,\,\int\limits_A
|\nabla u(x)|^n\,\,dm(x)\,,
\end{equation}
где $W_0(E)=W_0\left(A,\,C\right)$ -- семейство неотрицательных
непрерывных функций $u:A\rightarrow{\Bbb R}$ с компактным носителем
в $A,$ таких что $u(x)\ge 1$ при $x\in C$ и $u\in ACL.$
В формуле (\ref{eq1.1AB}), как обычно, $|\nabla
u|={\left(\sum\limits_{i=1}^n\,{\left(\partial_i u\right)}^2
\right)}^{1/2}.$
Говорят, что компакт $C$ в ${\Bbb R}^n,$ $n\ge 2,$ имеет {\it
нулевую ёмкость,} пишут ${\rm cap}\,C=0,$ если существует
ограниченное открытое множество $A,$ такое что $C\subset A$ и
ёмкость конденсатора $E=(A, C)$ удовлетворяет условию: ${\rm
cap}\,(A, C)=0.$ Известно, что в последнем случае и для любого
другого ограниченного открытого множества $A$ в ${\Bbb R}^n,$
содержащего $C,$ также будет выполнено условие: ${\rm cap}\,(A,
C)=0$ (см., напр., \cite[лемма~3.4, гл.~II]{Re}). Если хотя бы для
одного ограниченного открытого множества $A,$ такого что $C\subset
A,$ имеет место неравенство ${\rm cap}\,(A, C)>0,$ мы полагаем:
${\rm cap}\,C>0$ и в этом случае говорим, что множество $C$ имеет
{\it положительную ёмкость}.

\medskip\, Для понимания природы множеств ёмкости нуль полезно привести
следующее замечание (см. \cite[следствия~1-2, $\S\,3$]{Re$_1*$}, см.
также \cite[лемма~2.13]{MRV$_2$} и \cite[теорему~IV.4]{HW}). Пусть
$F\subset {\Bbb R}^n$ -- компактное множество нулевой ёмкости,
тогда: при каждом $\alpha
> 0$ его $\alpha$-мерная мера Хаусдорфа
${\mathcal H}^{\alpha}(F)$ равна нулю;  имеют место условия:
$m(F)=0$ и ${\rm Int}\,F = \varnothing;$ 3) множество $D\setminus F$
является областью (теорема Менгера--Урысона).

\medskip{} Пусть $\left(X,\,d\right)$ и
$\left(X^{\,{\prime}},{d}^{\,{\prime}}\right)$ -- метрические
пространства с расстояниями  $d$  и ${d}^{\,{\prime}}$
соответственно. Семейство $\frak{F}$ непрерывных отображений
$f:X\rightarrow {X}^{\,\,\prime}$ называется {\it нормальным}, если
из любой последовательности отображений $f_{m} \in \frak{F}$ можно
выделить подпоследовательность $f_{m_{k}}$, которая сходится
локально равномерно в $X$ к непрерывной функции
$f:\,X\,\rightarrow\, X^{\,\prime}.$

\medskip
Введенное  понятие  очень  тесно  связано  со  следующим. Семейство
$\frak{F}$ отображений $f:X\rightarrow {X}^{\,\prime}$ называется
{\it равностепенно непрерывным в точке} $x_0 \in X,$ если для любого
$\varepsilon>0$ найдётся $\delta>0$ такое, что ${d}^{\,\prime}
\left(f(x),f(x_0)\right)<\varepsilon$ для всех $x$ таких, что
$d(x,x_0)<\delta$ и для всех $f\in \frak{F}.$ Говорят, что
$\frak{F}$ {\it равностепенно непрерывно}, если $\frak{F}$
равностепенно непрерывно в каждой  точке $x_0\in X.$ Согласно одной
из версий теоремы Арцела-Асколи (см., напр., \cite[пункт~20.4]{Va}),
если $\left(X,\,d\right)$ -- сепарабельное метрическое пространство,
а $\left(X^{\,\prime},\, d^{\,\prime}\right)$ -- компактное
метрическое пространство, то семейство $\frak{F}$ отображений
$f:X\rightarrow {X}^{\,\prime}$ нормально тогда  и только тогда,
когда  $\frak{F}$ равностепенно непрерывно.

\medskip
Отметим, что всюду далее, если не оговорено противное, $(X, d)=(D,
|\cdot|),$ где $D$ -- область в ${\Bbb R}^n,$ а $|\cdot|$ --
евклидова метрика, $|x-y|=\sqrt{\sum\limits_{i=1}^n (y_i-x_i)^2},$
$x=(x_1,\ldots,x_n),$ $y=(y_1,\ldots,y_n);$ $\left(X^{\,\prime},\,
d^{\,\prime}\right)=\left(\overline{{\Bbb R}^n},\, h\right),$ где
$\overline{{\Bbb R}^n}={\Bbb R}^n\cup\{\infty\},$ $h$ -- хордальная
метрика,
$$h(x,\infty)=\frac{1}{\sqrt{1+{|x|}^2}}\,,\quad h(x,y)=\frac{|x-y|}{\sqrt{1+{|x|}^2} \sqrt{1+{|y|}^2}}\,,  x\ne
\infty\ne y\,.$$

\medskip{} Переходим к формулировке основных результатов настоящей работы.

\medskip
Для заданных компактного множества $E\subset{\Bbb R}^n,$ неубывающей
функции $\varphi:[0,\infty)\rightarrow[0,\infty),$ измеримой по
Лебегу функции $Q:D\rightarrow [1, \infty]$ и числа $N\in {\Bbb N}$
обозначим символом $\frak{R}_{\varphi, Q, N, E}$ семейство всех
открытых дискретных  отображений $f:D\rightarrow {\Bbb R}^n\setminus
E$ класса $W^{1,\varphi}_{loc},$ имеющих конечное искажение, таких
что $N(f , D)\le N$ и $K_O^{n-1}(x, f)\le Q(x)$ почти всюду.
Справедлива следующая

\medskip
\begin{theorem}\label{th1}
{\sl\, Пусть $n\ge 3,$ тогда семейство отображений
$\frak{R}_{\varphi, Q, N, E}$ является равностепенно непрерывным в
некоторой фиксированной точке $x_0\in D,$ если ${\rm cap}\,E>0,$
$Q\in L_{loc}^1(D),$
\begin{equation}\label{eqOS3.0a}
\int\limits_{1}^{\infty}\left[\frac{t}{\varphi(t)}\right]^
{\frac{1}{n-2}}dt<\infty
\end{equation}
и, кроме того, при некотором $\varepsilon_0>0,$ $\varepsilon_0<{\rm
dist}(x_0, \partial D),$ выполнено следующее условие расходимости
интеграла:
\begin{equation}\label{eq9}
\int\limits_{0}^{\varepsilon_0}
\frac{dt}{tq_{x_0}^{\,\frac{1}{n-1}}(t)}=\infty,
\end{equation}
где, как обычно,
\begin{equation}\label{eq9b}
q_{x_0}(r):=\frac{1}{\omega_{n-1}r^{n-1}}\int\limits_{|x-x_0|=r}Q(x)\,d{\mathcal
H}^{n-1}
\end{equation}
-- среднее интегральное значение функции $Q$ над сферой $S(x_0, r).$
В частности, заключение теоремы \ref{th1} является верным, если
$q_{x_0}(r)=\,O\left({\left[ \log{\frac{1}{r}}\right]}^{n-1}\right)$
при $r\rightarrow 0.$}
\end{theorem}

\medskip
Из приведённого выше критерия Арцела-Асколи вытекает следующее

\medskip
\begin{corollary}\label{cor2}{\sl\,
В условиях теоремы \ref{th1} семейство отображений
$\frak{R}_{\varphi, Q, N, E}$ является нормальным семейством
отображений, как только условие (\ref{eq9}) выполнено в каждой точке
$x_0$ области $D.$}
\end{corollary}

\medskip
Сформулируем ещё один важнейший результат работы.

\medskip
Будем говорить, что локально интегрируемая функция
${\varphi}:D\rightarrow{\Bbb R}$ имеет {\it конечное среднее
колебание} в точке $x_0\in D$, пишем $\varphi\in FMO(x_0),$ если
%
%
%
%
$\limsup\limits_{\varepsilon\rightarrow
0}\frac{1}{\Omega_n\varepsilon^n}\int\limits_{B( x_0,\,\varepsilon)}
|{\varphi}(x)-\overline{{\varphi}}_{\varepsilon}|\, dm(x)<\infty,$
%
%
где
$\overline{{\varphi}}_{\varepsilon}=\frac{1}
{\Omega_n\varepsilon^n}\int\limits_{B(x_0,\,\varepsilon)}
{\varphi}(x)\, dm(x).$
\medskip
Заметим, что, как известно, $\Omega_n\varepsilon^n=m(B(x_0,
\varepsilon)).$ Имеет место следующая

\medskip
\begin{theorem}\label{th2}
{\sl\, При $n\ge 3$ семейство отображений $\frak{R}_{\varphi, Q, N,
E}$ является равностепенно непрерывным в точке $x_0\in D,$ если
${\rm cap\,}E>0,$ выполнено условие (\ref{eqOS3.0a}) и, кроме того,
$Q\in FMO(x_0).$}
\end{theorem}

\medskip
Из теоремы \ref{th2} на основании приведённого выше критерия
Арцела-Асколи вытекает следующее

\medskip
\begin{corollary}\label{cor3}{\sl\,
В условиях теоремы \ref{th2} семейство отображений
$\frak{R}_{\varphi, Q, N, E}$ является нормальным семейством
отображений, как только условие $Q\in FMO(x_0)$ выполнено в каждой
точке $x_0$ области $D.$}
\end{corollary}

\section{Предварительные сведения}\label{sect2}

\medskip{} В данном параграфе обсуждаются различные не
связанные между собой вопросы, каждый из которых является
вспомогательным элементом при доказательстве основных результатов
работы. Прежде всего, напомним некоторые определения, связанные с
понятием поверхности, интеграла по поверхности, а также модулей
семейств кривых и поверхностей.

\medskip Пусть $\omega$ -- открытое множество в $\overline{{\Bbb
R}^k}:={\Bbb R}^k\cup\{\infty\},$ $k=1,\ldots,n-1.$ Непрерывное
отображение $S:\omega\rightarrow{\Bbb R}^n$ будем называть {\it
$k$-мерной поверхностью} $S$ в ${\Bbb R}^n.$ Число прообразов
$$N(y, S)={\rm card}\,S^{-1}(y)={\rm card}\,\{x\in\omega:S(x)=y\},\
y\in{\Bbb R}^n$$ будем называть {\it функцией кратности} поверхности
$S.$ Другими словами, $N(y, S)$ -- кратность накрытия точки $y$
поверхностью $S.$ Пусть $\rho:{\Bbb R}^n\rightarrow\overline{{\Bbb
R}^+}$ -- борелевская функция, в таком случае интеграл от функции
$\rho$ по поверхности $S$ определяется равенством:  $$\int\limits_S
\rho\,d{\mathcal{A}}:= \int\limits_{{\Bbb R}^n}\rho(y)\,N(y,
S)\,d{\mathcal H}^ky\,.$$
Пусть $\Gamma$ -- семейство $k$-мерных поверхностей $S.$ Борелевскую
функцию $\rho:{\Bbb R}^n\rightarrow\overline{{\Bbb R}^+}$ будем
называть {\it допустимой} для семейства $\Gamma,$ сокр. $\rho\in{\rm
adm}\,\Gamma,$ если
\begin{equation}\label{eq8.2.6}\int\limits_S\rho^k\,d{\mathcal{A}}\ge 1\end{equation}
для каждой поверхности $S\in\Gamma.$
Для заданного числа $p\in(0,\infty)$  {\it $p$-модулем} семейства
$\Gamma$ назовём величину
$$M_p(\Gamma)=\inf_{\rho\in{\rm adm}\,\Gamma} \int\limits_{{\Bbb
R}^n}\rho^p(x)\,dm(x)\,.$$ Мы также полагаем
$M(\Gamma)=M_n(\Gamma),$
а величину $M(\Gamma)$ в этом случае называем {\it модулем
семейства} $\Gamma.$ Заметим, что модуль семейств поверхностей,
определённый таким образом, представляет собой внешнюю меру в
пространстве всех $k$-мерных поверхностей (см. \cite{Fu}).

\medskip{} Пусть $p\ge 1.$ Говорят, что некоторое свойство $P$ выполнено для {\it $p$-почти всех
поверхностей} области $D,$ если оно имеет место для всех
поверхностей, лежащих в $D,$ кроме, быть может, некоторого их
подсемейства, $p$-модуль которого равен нулю. (Как правило, если
речь идёт о конформном модуле, говорят, что указанное свойство
выполнено для {\it почти всех поверхностей} области $D,$ опуская
приставку $"$$n$$"$ в выражении $"$$n$-почти всех$"$). В частности,
говорят, что некоторое свойство выполнено для {\it $p$-почти всех
кривых} области $D$, если оно имеет место для всех кривых, лежащих в
$D$, кроме, быть может, некоторого их подсемейства, $p$-модуль
которого равен нулю.

\medskip Будем говорить, что измеримая по Лебегу функция $\rho:{\Bbb
R}^n\rightarrow\overline{{\Bbb R}^+}$ {\it $p$-обобщённо допустима}
для семейства $\Gamma$ $k$-мерных поверхностей $S$ в ${\Bbb R}^n,$
сокр. $\rho\in{\rm ext}_p\,{\rm adm}\,\Gamma,$ если соотношение
(\ref{eq8.2.6}) выполнено для $p$-почти всех поверхностей $S$
семейства $\Gamma.$ {\it Обобщённый $p$-модуль} $\overline
M_p(\Gamma)$ семейства $\Gamma$ определяется равенством
$$\overline M_p(\Gamma)= \inf\int\limits_{{\Bbb
R}^n}\rho^p(x)\,dm(x)\,,$$
где точная нижняя грань берётся по всем функциям $\rho\in{\rm
ext}_p\,{\rm adm}\,\Gamma.$ В случае $p=n$ мы используем обозначения
$\overline M(\Gamma)$ и $\rho\in{\rm ext}\,{\rm adm}\,\Gamma,$
соответственно. Очевидно, что при каждом $p\in(0,\infty),$
$k=1,\ldots,n-1,$ и каждого семейства $k$-мерных поверхностей
$\Gamma$ в ${\Bbb R}^n,$ выполнено равенство $\overline
M_p(\Gamma)=M_p(\Gamma).$

\medskip{} Следующий класс отображений представляет собой обобщение
квазиконформных отображений в смысле кольцевого определения по
Герингу (\cite{Ge$_3$}) и отдельно исследуется различными авторами
(см., напр., \cite[глава~9]{MRSY}). Пусть $D$ и $D^{\,\prime}$ --
заданные области в $\overline{{\Bbb R}^n},$ $n\ge 2,$
$x_0\in\overline{D}\setminus\{\infty\}$ и $Q:D\rightarrow(0,\infty)$
-- измеримая по Лебегу функция. Будем говорить, что $f:D\rightarrow
D^{\,\prime}$ -- {\it нижнее $Q$-отображение в точке} $x_0,$ как
только
\begin{equation}\label{eq1A}
M(f(\Sigma_{\varepsilon}))\ge \inf\limits_{\rho\in{\rm
ext\,adm}\,\Sigma_{\varepsilon}}\int\limits_{D\cap
R_{\varepsilon}}\frac{\rho^n(x)}{Q(x)}\,dm(x)
\end{equation}
для каждого кольца $A(\varepsilon, \varepsilon_0, x_0),$
$\varepsilon_0\in(0,d_0),$ $d_0=\sup\limits_{z\in D}|z-z_0|,$
где $\Sigma_{\varepsilon}$ обозначает семейство всех пересечений
сфер $S(x_0, r)$ с областью $D,$ $r\in (0, \varepsilon_0).$ Мы не
будем пока что останавливаться на примерах таких отображений,
однако, заметим, что, как показывает один из основных результатов
настоящей работы, указанные примеры несложно указать (см. теоремы
\ref{thOS4.1} и \ref{thOS4.2}).

\medskip{}
Отметим, что выражения $"$почти всех кривых$"$ и $"$почти всех
по\-вер\-х\-но\-с\-тей$"$ в отдельных случаях могут иметь две
различные интерпретации. Например, в определении понятия $ACL$
выражение $"$почти всех отрезков$"$ необходимо понимать относительно
меры их проекции на соответствующую гиперплоскость. В то же время,
выражение $"$почти всех$"$ можно также интерпретировать в смысле
конформного модуля. Следующее утверждение вносит некоторую ясность
между указанными интерпретациями (см. по этому поводу также
\cite[лемма~9.1]{MRSY}).

\begin{lemma}\label{lem8.2.11}{\sl\, Пусть $x_0\in D.$ Если некоторое
свойство $P$ имеет место для почти всех сфер $D(x_0, r):=S(x_0,
r)\cap D,$ где $"$почти всех$"$ понимается в смысле модуля семейств
поверхностей, то $P$ также имеет место для почти всех сфер $D(x_0,
r)$ относительно линейной меры Лебега по параметру $r\in {\Bbb R }.$
Обратно, пусть $P$ имеет место для почти всех сфер $D(x_0,
r):=S(x_0, r)\cap D$ относительно линейной меры Лебега по $r\in
{\Bbb R},$ тогда $P$ также имеет место для почти всех поверхностей
$D(x_0, r):=S(x_0, r)\cap D$ в смысле модуля семейств
поверхностей.}\end{lemma}

\begin{proof} {\it Необходимость.} Пусть некоторое свойство
$P$ имеет место для почти всех сфер $D(x_0, r):=S(x_0, r)\cap D,$
где $"$почти всех$"$ понимается в смысле модуля семейств
поверхностей. Покажем, что $P$ также имеет место для почти всех сфер
$D(x_0, r)$ по отношению к параметру $r\in {\Bbb R}.$

Достаточно рассмотреть случай, когда область $D$ ограничена.
Предположим, что заключение леммы не является верным. Тогда найдётся
семейство $\Gamma$ сфер $D(x_0, r),$ для которого свойство $P$
выполнено в смысле почти всех поверхностей относительно модуля,
однако, нарушается для некоторого множества индексов $r\in {\Bbb R}$
положительной меры.

Ввиду регулярности меры Лебега $m_{1}$ найдётся борелевское
множество $B\subset {\Bbb R}$ такое что $m_{1}(B)>0$ и свойство $P$
нарушается для почти всех $r\in B.$ Пусть $\rho:{\Bbb
R}^n\rightarrow[0,\infty]$ -- допустимая функция для семейства
$\Gamma.$ Учитывая, что $B$ -- борелево, мы можем считать, что
$\rho\equiv 0$ вне $E=\{x\in D:\exists\,\, r\in B: |x-x_0|=r\},$
поскольку, в этом случае, множество $E,$ очевидно, борелево. По
неравенству Гёльдера
$$\int\limits_{E}\rho^{n-1}(x)\ dm(x)\ \le\
\left(\int\limits_{E}\rho^n(x)\
dm(x)\right)^{\frac{n-1}{n}}\left(\int\limits_{E}\
dm(x)\right)^{\frac{1}{n}}$$
и следовательно, ввиду теоремы Фубини (см. \cite[Теорема~8.1,
гл.~III]{Sa}),
$$\int\limits_{{\Bbb R}^n}\rho^n(x)\ dm(x)\ \ge\ \frac{\left(\int\limits_{E}\rho^{n-1}(x)\
dm(x)\right)^{\frac{n}{n-1}}}{\left(\int\limits_{E}\
dm(x)\right)^{\frac{1}{n-1}}}\ \ge\
\frac{(m_{1}(B))^{\frac{n}{n-1}}}{c}$$ для некоторого $c>0,$ т.е.,
$M(\Gamma)>0,$  что противоречит предположению леммы. Первая часть
леммы \ref{lem8.2.11} доказана.

\medskip
{\it Достаточность.} Пусть $P$ имеет место для почти всех $r$
относительно меры Лебега и всех соответствующих этим $r$ сфер
$D(x_0, r),$ $r\in {\Bbb R}.$ Покажем, что $P$ также выполняется для
почти всех поверхностей $D(x_0, r):=S(x_0, r)\cap D$ в смысле модуля
семейств поверхностей.

Обозначим через $\Gamma_0$ семейство всех пересечений
$D_r:=D(x_0,r)$ сфер $S(x_0, r)$ с областью $D,$ для которых $P$ не
имеет места. Пусть $R$ обозначает множество всех $r\in {\Bbb R}$
таких, что $D_r\in\Gamma_0.$ Если $m_1(R)=0,$ то по теореме Фубини
мы получаем, что $m(E)=0,$ где $E=\{x\in D: |x-x_0|=r\in R\}.$
Рассмотрим функцию $\rho_1:{\Bbb R}^n\rightarrow [0, \infty],$
определённую символом $\infty$ при $x\in E,$ и доопределённая нулём
в остальных точках. Отметим, что найдётся борелева функция
$\rho_2:{\Bbb R}^n\rightarrow [0, \infty],$ совпадающая почти всюду
с $\rho_1$ (см. \cite[раздел~2.3.5]{Fe}). Таким образом,
$M(\Gamma_0)\le \int\limits_E \rho_2^n dm(x)=\int\limits_E \rho_1^n
dm(x)=0,$ следовательно, $M(\Gamma_0)=0.$ Лемма \ref{lem8.2.11}
полностью доказана.
\end{proof}

\medskip{} Следующее утверждение не имеет прямого отношения к исследованию отображений,
однако, как будет видно далее, является весьма полезным (см.
\cite[лемма~9.2]{MRSY}).

\medskip
\begin{proposition}\label{Salnizh1} {\sl Пусть $(X, \mu)$ -- измеримое пространство с
конечной мерой $\mu,$ $q\in(1,\infty),$ и пусть
$\varphi:X\to(0,\infty)$ -- измеримая функция. Полагаем
\begin{equation}\label{Sal_eq2.1.7}
I(\varphi, q)=\inf\limits_{\alpha}
\int\limits_{X}\varphi\,\alpha^q\,d\mu\,,
\end{equation}
где инфимум берется по всем измеримым функциям
$\alpha:X\rightarrow[0,\infty]$ таким, что
%
$\int\limits_{X}\alpha\,d\mu=1.$ 
Тогда
%
$I(\varphi, q)=\left[\int\limits_{X}\varphi^{-\lambda}\,d\mu\right]
^{-\frac{1}{\lambda}},$ 
где
$\lambda=\frac{q^{\,\prime}}{q},$ $\frac{1}{q}+\frac{1}{q^{\,\prime}}=1,$ 
т.е. $\lambda=1/(q-1)\in(0,\infty).$ Точная нижняя грань в
(\ref{Sal_eq2.1.7}) достигается на функции
$$
\rho=\left(\int\limits_X
\varphi^\frac{1}{1-\alpha}\;d\mu\right)^{-1}\varphi^\frac{1}{1-\alpha}\,.
$$}
\end{proposition}

\medskip{} Обратимся теперь к изучению соотношений вида (\ref{eq1A}).
Отметим, что условие (\ref{eq1A}) представляет собой бесконечную
серию неравенств по $\rho$ и потому в большинстве случаев может
оказаться весьма затруднительным для проверки. Следующее утверждение
облегчает решение указанной проблемы ввиду возможности проверки
всего одного неравенства по отношению к фиксированному семейству
$\Sigma_{\varepsilon},$ а не бесконечного их ряда. Имеет место
следующее утверждение, являющееся обобщением теоремы 9.2 в
\cite{MRSY}.

\medskip
\begin{lemma}\label{lem4}{\sl\,
Пусть $D,$  $D^{\,\prime}\subset\overline{{\Bbb R}^n},$
$x_0\in\overline{D}\setminus\{\infty\}$ и $Q:D\rightarrow(0,\infty)$
-- измеримая по Лебегу функция. Отображение $f:D\rightarrow
D^{\,\prime}$ является нижним $Q$-отображением в точке $x_0$ тогда и
только тогда, когда
\begin{equation}\label{eq15}
M(f(\Sigma_{\varepsilon}))\ge\int\limits_{\varepsilon}^{\varepsilon_0}
\frac{dr}{||\,Q||_{n-1}(r)}\quad\forall\
\varepsilon\in(0,\varepsilon_0)\,,\ \varepsilon_0\in(0,d_0)\,,
\end{equation}
где, как и выше, $\Sigma_{\varepsilon}$ обозначает семейство всех
пересечений сфер $S(x_0, r)$ с областью $D,$ $r\in (0,
\varepsilon_0),$
$$
\Vert
Q\Vert_{n-1}(r)=\left(\int\limits_{D(x_0,r)}Q^{n-1}(x)\,d{\mathcal{A}}\right)^{\frac{1}{n-1}}$$
-- $L_{n-1}$-норма функции $Q$ над сферой $D(x_0,r)=\{x\in D:
|x-x_0|=r\}=D\cap S(x_0,r)$.}
\end{lemma}

\medskip
\begin{proof} Не ограничивая общности рассуждений, можно считать, что
$\Vert Q\Vert_{n-1}(r)<\infty$ при почти всех $r\in (\varepsilon,
\varepsilon_0).$ Рассмотрим произвольное кольцо $A(\varepsilon,
\varepsilon_0, x_0),$ где $\varepsilon_0\in(0,d_0),$
$d_0=\sup\limits_{x\in D}|x-x_0|,$ а $\Sigma_{\varepsilon}$
обозначает семейство всех пересечений сфер $S(x_0, r)$ с областью
$D,$ $r\in (0, \varepsilon_0).$ Пусть $\rho\in {\rm
ext\,adm}\,\Sigma_{\varepsilon},$ тогда по лемме \ref{lem8.2.11}
$$A_{\rho}(r)=\int\limits_{D(x_0,r)}\rho^{n-1}(x)\ d{\mathcal A}\ne 0$$
при почти всех $r\in (0, \varepsilon_0)$ и $A_{\rho}(r)$ --
измеримая по Лебегу функция относительно параметра $r$ ввиду теоремы
Фубини. Пусть $A_{\Sigma_\varepsilon}$ обозначает класс всех
измеримых по Лебегу функций $\rho:{\Bbb
R}^n\rightarrow\overline{{\Bbb R}^+},$ удовлетворяющих условию
$\int\limits_{D(x_0, r)}\rho^{n-1}\,d{\mathcal A}=1$ для почти всех
$r\in (0, \varepsilon_0).$ Поскольку $A_{\Sigma_\varepsilon}\subset
{\rm ext\,adm}\,\Sigma_{\varepsilon},$  мы получим, что
\begin{equation}\label{eq1}
\inf\limits_{\rho\in{\rm
ext\,adm}\,\Sigma_{\varepsilon}}\int\limits_{D\cap
R_{\varepsilon}}\frac{\rho^n(x)}{Q(x)}\,dm(x)\le
\inf\limits_{\rho\in A_{\Sigma_\varepsilon}}\int\limits_{D\cap
R_{\varepsilon}}\frac{\rho^n(x)}{Q(x)}\,dm(x)\,.
\end{equation}
С другой стороны, для заданной функции $\rho\in {\rm
ext\,adm}\,\Sigma_{\varepsilon}$ мы получим, что
$$\inf\limits_{\rho\in A_{\Sigma_\varepsilon}}\int\limits_{D\cap
R_{\varepsilon}}\frac{\rho^n(x)}{Q(x)}\,dm(x)\le
\int\limits_{\varepsilon}^{\varepsilon_0}\frac{1}{\left(A_{\rho}(r)\right)^
{\frac{n}{n-1}}}\int\limits_{D(x_0, r)}\frac{\rho^n(x)}{Q(x)}\,
d{\mathcal A}dr\le$$
\begin{equation}\label{eq2}
\le \int\limits_{D\cap R_{\varepsilon}}\frac{\rho^n(x)}{Q(x)}\,dm(x)
\end{equation}
поскольку $A_{\rho}(r)\ge 1$ для почти всех $r\in (0,
\varepsilon_0)$ ввиду леммы \ref{lem8.2.11}. Из (\ref{eq1}) и
(\ref{eq2}) мы будем иметь, что
\begin{equation}\label{eq12} \inf\limits_{\rho\in{\rm ext}\,{\rm
adm}\,\Sigma_{\varepsilon}}\int\limits_{D\cap
R_{\varepsilon}}\frac{\rho^n(x)}{Q(x)}\,dm(x)=\inf\limits_{\rho\in
A_{\Sigma_\varepsilon}}\int\limits_{D\cap
R_{\varepsilon}}\frac{\rho^n(x)}{Q(x)}\,dm(x)\,.
\end{equation}
Теперь покажем, что
\begin{equation}\label{eq9A}
\inf\limits_{\rho\in{\rm ext}\,{\rm
adm}\Sigma_{\varepsilon}}\int\limits_{D\cap
R_{\varepsilon}}\frac{\rho^n(x)}{Q(x)}\ dm(x)\ =\
\int\limits_{\varepsilon}^{\varepsilon_0}\left(\inf\limits_{\alpha\in
I(r)}\int\limits_{D(x_0,r)}\frac{\alpha^q(x)}{Q(x)}\ d{\mathcal
H}^{n-1}\right)dr\,,
\end{equation}
где $q=n/(n-1)>1$ и $I(r)$ обозначает множество всех измеримых
функций $\alpha$ на сфере $D(x_0,r)=S(x_0,r)\cap D$ таких что
$$\int\limits_{D(x_0,r)}\alpha(x) d{\mathcal H}^{n-1}=1\,.$$

\medskip
Прежде всего, обозначая $\psi(r):=\inf\limits_{\alpha\in
I(r)}\int\limits_{D(x_0,r)}\frac{\alpha^q(x)}{Q(x)}\ d{\mathcal
H}^{n-1},$ ввиду леммы \ref{lem1} мы будем иметь, что
\begin{equation}\label{eq14}
\psi(r)=(\Vert Q\Vert_{n-1}(r))^{-1}=\left(\int\limits_{D(x_0,r)}
Q^{n-1}(x)\,d{\mathcal H}^{n-1}x\right)^{-\frac{1}{n-1}}\,.
\end{equation}
Следовательно, ввиду теоремы Фубини функция $\psi(r)$ измерима по
$r$ (см., напр., \cite[теорема~8.1, гл.~III]{Sa}), так что интеграл
в правой части (\ref{eq9A}) определён корректно.

\medskip
Пусть $\rho\in A_{\Sigma_\varepsilon},$ тогда функция
$\rho_r(x):=\rho|_{S(x_0, r)}$ измерима по отношению к хаусдорфовой
мере ${\mathcal H}^{n-1}$ для п.в. $r\in (\varepsilon,
\varepsilon_0)$ ввиду теоремы Фубини (см. \cite[теорема~8.1,
гл.~III]{Sa}). Следовательно,
$$\int\limits_{D\cap R_{\varepsilon}}\frac{\rho^n(x)}{Q(x)} dm(x)=
\int\limits_{\varepsilon}^{\varepsilon_0}\int\limits_{D(x_0,
r)}\frac{\rho_r^n(x)}{Q(x)}\ d{\mathcal H}^{n-1}dr\ge$$
%
$$\ge\int\limits_{\varepsilon}^{\varepsilon_0}\left(\inf\limits_{\alpha\in
I(r)}\int\limits_{D(x_0,r)}\frac{\alpha^q(x)}{Q(x)}\ d{\mathcal
H}^{n-1}\right)dr\,.$$
%
Переходя к $\inf$ по всем $\rho\in A_{\Sigma_\varepsilon},$ ввиду
(\ref{eq12}) мы получим, что
\begin{equation}\label{eq11}
\inf\limits_{\rho\in {\rm ext}\,{\rm
adm}\Sigma_{\varepsilon}}\int\limits_{D\cap
R_{\varepsilon}}\frac{\rho^n(x)}{Q(x)}\ dm(x)\ge
\int\limits_{\varepsilon}^{\varepsilon_0}\left(\inf\limits_{\alpha\in
I(r)}\int\limits_{D(x_0,r)}\frac{\alpha^q(x)}{Q(x)}\ d{\mathcal
H}^{n-1}\right)dr\,.
\end{equation}
Докажем теперь обратное неравенство. Ввиду предложения
\ref{Salnizh1} точная нижняя грань функции
$\psi(r)=\int\limits_{D(x_0,r)}\frac{\alpha^q(x)}{Q(x)}\ d{\mathcal
H}^{n-1}$ в (\ref{eq11}) достигается на функции
$$\alpha_0(x):=
\frac{Q^{n-1}(x)}{\int\limits_{D(x_0,r)} Q^{n-1}(x)\,d{\mathcal
H}^{n-1}x}\,.$$
Ввиду сделанного выше предположения $\Vert Q\Vert_{n-1}(r)<\infty$
при почти всех $r\in (\varepsilon, \varepsilon_0).$ Следовательно,
$\alpha_0^{\frac{1}{n-1}}\in A_{\Sigma_\varepsilon}$ и, значит,
$$\inf\limits_{\rho\in {\rm ext}\,{\rm
adm}\Sigma_{\varepsilon}}\int\limits_{D\cap
R_{\varepsilon}}\frac{\rho^n(x)}{Q(x)}\ dm(x)\le$$
\begin{equation}\label{eq13}
\le \int\limits_{D\cap
R_{\varepsilon}}\frac{\alpha_0^{\frac{n}{n-1}}(x)}{Q(x)}\ dm(x)=
\int\limits_{\varepsilon}^{\varepsilon_0}\left(\inf\limits_{\alpha\in
I(r)}\int\limits_{D(x_0,r)}\frac{\alpha^q(x)}{Q(x)}\ d{\mathcal
H}^{n-1}\right)dr\,.
\end{equation}
Из неравенств (\ref{eq11}) и (\ref{eq13}) вытекает соотношение
(\ref{eq9A}). С другой стороны, (\ref{eq15}) вытекает из
(\ref{eq9A}) и (\ref{eq14}).
\end{proof}

\medskip{}
Следующее утверждение имеет важное значение для доказательства
многих результатов настоящей работы (см. \cite[предложение~10.2,
гл.~II]{Ri}).

\medskip
\begin{proposition}\label{pr1*!}{\,\sl Пусть $E=(A,\,C)$ --
произвольный конденсатор в ${\Bbb R}^n$ и пусть $\Gamma_E$ --
семейство всех кривых вида $\gamma:[a,\,b)\rightarrow A$ таких, что
$\gamma(a)\in C$ и $|\gamma|\cap\left(A\setminus
F\right)\ne\varnothing$ для произвольного компакта $F\subset A.$
Тогда
%
$${\rm cap}\,E=M(\Gamma_E)\,.$$
%
}
\end{proposition}

\medskip{}
В дальнейшем всюду символом $\Gamma(E,F,D)$ мы обозначаем семейство
всех кривых $\gamma:[a,b]\rightarrow\overline{{\Bbb R}^n},$ которые
соединяют $E$ и $F$ в $D,$ т.е. $\gamma(a)\in E,\,\gamma(b)\in F$ и
$\gamma(t)\in D$ при $t\in(a,\,b).$ Для доказательства основных
результатов работы также существенно используются так называемые
кольцевые $Q$-отображения, определение которых приведено ниже (см.,
напр., \cite[гл.~7]{MRSY}, см. также \cite{BGMV} и \cite{GRSY}).
Говорят, что $f:D\rightarrow \overline{{\Bbb R}^n}$ является {\it
кольцевым $Q$-отоб\-ра\-же\-нием в точке $x_0\,\in\,D,$} если
соотношение
%
$$M\left(f\left(\Gamma\left(S_1,\,S_2,\,A\right)\right)\right)\ \le
\int\limits_{A} Q(x)\cdot \eta^n(|x-x_0|)\ dm(x)$$ 
выполнено для любого кольца $A=A(r_1,r_2, x_0),$\, $0<r_1<r_2<
r_0:={\rm dist\,}(x_0, \partial D),$ и для каждой измеримой функции
$\eta : (r_1,r_2)\rightarrow [0,\infty ]\,$ такой, что
%
%
%
$\int\limits_{r_1}^{r_2}\eta(r)dr\ge 1.$
Отметим, что кольцевые $Q$-отображения являются обобщением
квазиконформных отображений и отображений с ограниченным искажением
(см., напр., \cite{BGMV}, \cite{GRSY}, \cite{Pol}, \cite{Re$_1*$},
\cite{Re}, \cite{Ri}, \cite{Va} и \cite{Vu}). В частности, ввиду
неравенства Е.~Полецкого отображения с ограниченным искажением
являются кольцевыми $Q$-отоб\-ра\-же\-ни\-я\-ми с некоторой
ограниченной функцией $Q$ (см. \cite[теорема~1]{Pol}). Следующее
утверждение было доказано и опубликовано автором данной работы
несколько ранее и может быть найдено в \cite[теорема~1]{Sev}.

\medskip
\begin{lemma}\label{lem2} {\sl\, Пусть $Q:D\rightarrow [0, \infty]$ -- измеримая по Лебегу
функция, $Q\in L_{loc}^1(D).$ Открытое дискретное отображение
$f:D\rightarrow \overline{{\Bbb R}^n}$ является кольцевым
$Q$-отображением в точке $x_0\in D$ тогда и только тогда, когда для
произвольных $0<r_1<r_2< {\rm dist} \, (x_0,\partial D)$ и
произвольного конденсатора $E=\left(B(x_0, r_2), \overline{B(x_0,
r_1)}\right)$ ёмкость конденсатора $$f(E):=\left(f(B(x_0, r_2)),
f\left(\overline{B(x_0, r_1)}\right)\right)$$ удовлетворяет условию
%
$${\rm cap}\, f(E)\le\frac{\omega_{n-1}}{I^{n-1}}\,,$$
%
где $I=I(x_0,r_1,r_2)$ задаётся соотношением
%
$$I=I(x_0,r_1,r_2)=\int\limits_{r_1}^{r_2}\
\frac{dr}{rq_{x_0}^{\frac{1}{n-1}}(r)}\,.$$
}
\end{lemma}

\medskip{} Следующие важные сведения, касающиеся ёмкости пары
множеств относительно области, могут быть найдены в работе В.~Цимера
\cite{Zi}. Пусть $G$ -- ограниченная область в ${\Bbb R}^n$ и $C_{0}
, C_{1}$ -- непересекающиеся компактные множества, лежащие в
замыкании $G.$ Полагаем  $R=G \setminus (C_{0} \cup C_{1})$ и
$R^{\,*}=R \cup C_{0}\cup C_{1}.$ {\it Конформной ёмкостью пары
$C_{0}, C_{1}$ относительно замыкания $G$} называется величина
$$C[G, C_{0}, C_{1}] = \inf \int\limits_{R} \vert \nabla u
\vert^{n}\ dm(x),$$
где точная нижняя грань берётся по всем функциям $u,$ непрерывным в
$R^{\,*},$ $u\in ACL(R),$ таким что $u=1$ на $C_{1}$ и $u=0$ на
$C_{0}.$ Указанные функции будем называть {\it допустимыми} для
величины $C [G, C_{0}, C_{1}].$ Мы будем говорить, что  {\it
множество $\sigma \subset {\Bbb R}^n$ разделяет $C_{0}$ и $C_{1}$ в
$R^{\,*}$} если $\sigma \cap R$ замкнуто в $R$ и найдутся
непересекающиеся множества $A$ и $B,$ являющиеся открытыми в
$R^{\,*} \setminus \sigma,$ такие что $R^{\,*} \setminus \sigma =
A\cup B,$ $C_{0}\subset A$ и $C_{1} \subset B.$ Пусть $\Sigma$
обозначает класс всех множеств, разделяющих $C_{0}$ и $C_{1}$ в
$R^{\,*}.$ Для числа $n^{\prime} = n/(n-1)$ определим величину
\begin{equation}\label{eq13.4.12}
\widetilde{M_{n^{\prime}}}(\Sigma)=\inf\limits_{\rho\in
\widetilde{\rm adm} \Sigma} \int\limits_{{\Bbb
R}^n}\rho^{\,n^{\prime}}dm(x)
\end{equation}
где запись $\rho\in \widetilde{\rm adm}\,\Sigma$ означает, что
$\rho$ -- неотрицательная борелевская функция в ${\Bbb R}^n$ такая,
что
\begin{equation} \label{eq13.4.13}
\int\limits_{\sigma \cap R}\rho d{\mathcal H}^{n-1} \ge
1\quad\forall\, \sigma \in \Sigma\,. \end{equation}
Заметим, что согласно результата Цимера,
\begin{equation}\label{eq3}
\widetilde{M_{n^{\,\prime}}}(\Sigma)=C [G , C_{0} ,
C_{1}]^{\,-1/(n-1)}\,,
\end{equation}
см. \cite[теорема~3.13]{Zi}. Заметим также, что согласно результата
Хессе
\begin{equation}\label{eq4}
M(\Gamma(E, F, D))= C[D, E, F]\,,
\end{equation}
как только $(E \cap F)\cap
\partial D = \varnothing,$
см. \cite[теорема~5.5]{Hes}.

\section{Основная лемма}\label{sect3}

\medskip{} Дальнейшие исследования работы опираются на две важные взаимосвязи:
указывается взаимосвязь нижних $Q$-отоб\-ра\-же\-ний и кольцевых
$Q$-отоб\-ра\-же\-ний, а затем -- взаимосвязь классов
Орлича--Соболева с нижними $Q$-отоб\-ра\-же\-ни\-я\-ми. Подобные
взаимосвязи найдены нами, правда, только при дополнительном
предположении открытости и дискретности отображений рассматриваемых
классов. (Изучение более общего случая не относится к ближайшим
целям нашего исследования и, вероятно, требует подходов, существенно
отличных от модульной техники). Следующая лемма является наиболее
важным элементом дальнейшего изложения. Именно, здесь сформулирована
и доказана взаимосвязь нижних и кольцевых $Q$-отображений во
внутренних точках.

\medskip
\begin{lemma}\label{lem1} {\sl Пусть $x_0\in D$ и $Q:D\rightarrow
[0,\infty]$ -- локально интегрируемая в степени $n-1$ в $D$ функция.
Если $f:D\rightarrow \overline{{\Bbb R}^n}$ -- открытое дискретное
нижнее $Q$-отображение в точке $x_0,$ то $f$ является кольцевым
$Q^{\,*}$-отображением в этой же точке при
$Q^{\,*}=Q^{n-1}.$}\end{lemma}

\begin{proof} Пусть $x_0\in D,$ $0<r_1<r_2<{\rm dist\,}(x_0, \partial D).$
Без ограничения общности, мы можем считать, что $f(x_0)\ne \infty.$
Согласно лемме \ref{lem2} достаточно установить, что
$${\rm cap}\,
f(E)\le \frac{\omega_{n-1}}{I^{*\,n-1}}$$
где $E$ -- конденсатор вида $E=(B(x_0, r_2), \overline{B(x_0,
r_1)}),$ $\omega_{n-1}$ -- площадь единичной сферы в ${\Bbb R}^n,$
$q^{\,*}_{x_0}(r)$ -- среднее значение функции $Q^{n-1}(x)$ над
сферой $|x-x_0|=r$ и $I^{\,*}=I^{\,*}(x_0,
r_1,r_2)=\int\limits_{r_1}^{r_2}\
\frac{dr}{rq^{\,*\,\frac{1}{n-1}}_{x_0}(r)}.$ Зафиксируем
$\varepsilon\in (r_1, r_2)$ и рассмотрим шар $B(x_0, \varepsilon).$
Полагаем $C_0=\partial f(B(x_0, r_2)),$ $C_1=f(\overline{B(x_0,
r_1)}),$ $\sigma=\partial f(B(x_0, \varepsilon)).$ Поскольку
$\overline{B(x_0, r_2)}$ -- компакт в $D,$ найдётся шар $B(x_0, R)$
такой, что $\overline{f(B(x_0, r_2))}\subset B(x_0, R).$ Полагаем
$G:=B(x_0, R).$

\medskip
Поскольку $f$ -- непрерывно и открыто, $\overline{f(B(x_0, r_1))}$
-- компактное подмножество $f(B(x_0, \varepsilon))$ также, как
$\overline{f(B(x_0, \varepsilon))}$ -- компактное подмножество
$f(B(x_0, r_2)).$ В частности, $\overline{f(B(x_0, r_1))}\cap
\partial f(B(x_0, \varepsilon))=\varnothing.$ Пусть, как и выше, $R=G
\setminus (C_{0} \cup C_{1})$ и $R^{\,*} = R \cup C_{0}\cup C_{1},$
тогда $R^{\,*}:=G.$ Заметим, что $\sigma$ разделяет $C_0$ и $C_1$ в
$R^{\,*}=G.$ Действительно, множество $\sigma \cap R$ замкнуто в
$R,$ кроме того, пусть $A:=G\setminus \overline{f(B(x_0,
\varepsilon))}$ и $B= f(B(x_0, \varepsilon)),$ тогда $A$ и $B$
открыты в $G\setminus \sigma,$ $C_0\subset A,$ $C_1\subset B$ и
$G\setminus \sigma=A\cup B.$

\medskip
Пусть $\Sigma$ -- семейство всех множеств, отделяющих $C_0$ от $C_1$
в $G.$ Поскольку для открытых отображений $\partial f(O)\subset
f(\partial O),$ где $O$ -- компактная подобласть $D,$ мы получим:
$\partial f(B(x_0, r))\subset f(\partial B(x_0, r)),$ $r\in (0, {\rm
dist\,}(x_0, \partial D)).$

\medskip
Пусть $\rho^{n-1}\in \widetilde{{\rm
adm}}\bigcup\limits_{r=r_1}^{r_2}
\partial f(B(x_0, r))$ в смысле соотношения (\ref{eq13.4.13}), тогда также
$\rho\in {\rm adm}\bigcup\limits_{r=r_1}^{r_2}
\partial f(B(x_0, r))$ в смысле соотношения (\ref{eq8.2.6}).
Поскольку (ввиду открытости отображения $f$) имеет место включение
$\partial f(B(x_0, r))\subset f(S(x_0, r)),$ мы получим, что
$\rho\in {\rm adm}\bigcup\limits_{r=r_1}^{r_2} f(S(x_0, r))$ и,
следовательно, ввиду (\ref{eq13.4.12}) будем иметь
$$\widetilde{M_{n^{\prime}}}(\Sigma)\ge
\widetilde{M_{n^{\prime}}}\left(\bigcup\limits_{r=r_1}^{r_2}
\partial f(B(x_0, r))\right)\ge M\left(\bigcup\limits_{r=r_1}^{r_2}
\partial f(B(x_0, r))\right)\ge
$$
\begin{equation}\label{eq5}
\ge M\left(\bigcup\limits_{r=r_1}^{r_2} f(S(x_0, r))\right)\,.
\end{equation}
Однако, ввиду (\ref{eq3}) и (\ref{eq4}),
\begin{equation}\label{eq6}
\widetilde{M_{n^{\prime}}}(\Sigma)=\frac{1}{(M(\Gamma(C_0, C_1,
G)))^{1/(n-1)}}\,.
\end{equation}
Пусть $\Gamma_{f(E)}$ -- семейство всех кривых для конденсатора
$f(E)$ в обозначениях предложения \ref{pr1*!}. Пусть также
$\Gamma^{\,*}_{f(E)}$ обозначает семейство всех спрямляемых кривых
семейства $\Gamma_{f(E)},$ тогда заметим, что семейства
$\Gamma^{*}_{f(E)}$ и $\Gamma (C_0, C_1, G)$ имеют одинаковые
семейства допустимых метрик $\rho$ и, значит,
$$M(\Gamma_{f(E)})=M(\Gamma(C_0, C_1, G))\,.$$ Из (\ref{eq6}) и
предложения \ref{pr1*!} мы получим, что
\begin{equation}\label{eq7}
\widetilde{M}^{n-1}(\Sigma)=\frac{1}{{\rm cap\,}f(E)}\,.
\end{equation}
Окончательно, из (\ref{eq5}) и (\ref{eq7}) мы получаем неравенство
\begin{equation}\label{eq8}
{\rm cap\,}f(E) \le \frac{1}{M\left(\bigcup\limits_{r=r_1}^{r_2}
f(S(x_0, r))\right)^{n-1}}\,.
\end{equation}
По лемме \ref{lem4} и из (\ref{eq8}) мы получим, что
$$
{\rm cap\,}f(E) \le \frac{1}{\left(\int\limits_{r_1}^{r_2}
\frac{dr}{\Vert
\,Q\Vert_{n-1}(r)}\right)^{n-1}}=\frac{\omega_{n-1}}{I^{*\,n-1}}\,,
$$
что и доказывает утверждение леммы \ref{lem1}.
\end{proof}

\section{Взаимосвязь классов Орлича--Соболева с нижними
$Q$-отображениями}\label{sect4}

\medskip{} Результаты, сформулированные ниже, позволяют исследовать
классы Соболева и Орлича-Соболева, допускающие наличие точек
ветвления. В частности, ниже будет указана взаимосвязь нижних
$Q$-отображений с классами Соболева $W_{loc}^{1, 1}$ на плоскости.

\medskip{} Прежде всего, опишем взаимосвязь классов
Орлича--Соболева с нижними $Q$-отображениями при $n\ge 3.$

\medskip
Напомним, что отображение $f:X\rightarrow Y$ между пространствами с
мерами $(X, \Sigma, \mu)$ и $(Y, \Sigma^{\,\prime}, \mu^{\,\prime})$
обладает {\it $N$-свой\-с\-т\-вом} (Лузина), если из условия
$\mu(S)=0$ следует, что $\mu^{\,\prime}(f(S))=0.$ Следующее
вспомогательное утверждение получено в работе \cite{KRSS} (см.
теорема 1 и следствие 2).

\medskip
\begin{proposition}\label{pr1}
{\sl\, Пусть $D$ -- область в ${\Bbb R}^n,$ $n\ge 3,$
$\varphi:(0,\infty)\rightarrow (0,\infty)$ -- неубывающая функция,
удовлетворяющая условию (\ref{eqOS3.0a}). Тогда:

1) Если $f:D\rightarrow{\Bbb R}^n$ -- непрерывное открытое
отображение класса $W^{1,\varphi}_{loc}(D),$ то $f$ имеет почти
всюду полный дифференциал в $D;$

2) Любое непрерывное отображение $f\in W^{1,\varphi}_{loc}$ обладает
$N$-свойством относительно $(n-1)$-мерной меры Хаусдорфа, более
того, локально абсолютно непрерывно на почти всех сферах $S(x_0, r)$
с центром в заданной предписанной точке $x_0\in{\Bbb R}^n$. Кроме
того, на почти всех таких сферах $S(x_0, r)$ выполнено условие
${\mathcal H}^{n-1}(f(E))=0,$ как только $|\nabla f|=0$ на множестве
$E\subset S(x_0, r).$ (Здесь $"$почти всех$"$ понимается
относительно линейной меры Лебега по параметру $r$).}

\end{proposition}

Следующее утверждение может быть доказано по аналогии с
\cite[теорема~5]{KRSS}.

\medskip
\begin{theorem}{}\label{thOS4.1} {\sl Пусть $D$ -- область в ${\Bbb R}^n,$
$n\ge 3,$ $\varphi:(0,\infty)\rightarrow (0,\infty)$ -- неубывающая
функция, удовлетворяющая условию (\ref{eqOS3.0a}).
Если $n\ge 3,$ то каждое открытое дискретное отображение
$f:D\rightarrow {\Bbb R}^n$ с конечным искажением класса
$W^{1,\varphi}_{loc}$ такое, что $N(f, D)<\infty,$ является нижним
$Q$-отображением в каждой точке $x_0\in\overline{D}$ при $Q(x)=N(f,
D)\cdot K_O(x, f),$ где внешняя дилатация $K_O(x, f)$ отображения
$f$ в точке $x$ определена соотношением (\ref{eq0.1.1A}), а
кратность $N(f, D)$ определена соотношением (\ref{eq1.7A}).}
\end{theorem}

\medskip{} Пусть $D\subset {\Bbb C}.$ Для комплекснозначной
функции $f:D\rightarrow {\Bbb C},$ заданной в области $D\subset
{\Bbb C},$ имеющей частные производные по $x$ и $y$ при почти всех
$z=x + iy,$ полагаем $\overline{\partial} f= f_{\overline{z}} =
\left(f_x + if_y\right)/2$ и $\partial f = f_z = \left(f_x -
if_y\right)/2.$

\medskip
Отметим, что при $n=2$ связь классов Орлича--Соболева с нижними
$Q$-отоб\-ра\-же\-ни\-я\-ми значительно более проста, чем в
пространственном случае, как показывает теорема, приведённая ниже.

\medskip
\begin{theorem}\label{thOS4.2}{\sl\, Каждое открытое дискретное отображение $f:D\rightarrow
{\Bbb C}$ конечного искажения класса $W^{1, 1}_{\rm loc}$ такое, что
$N(f, D)<\infty,$ является нижним $Q$-отображением в произвольной
точке $z_0\in\overline{D}$ при $Q(z)=N(f, D)K_{\mu}(z),$ где
$K_{\mu}(z)$ определено соотношением
\begin{equation}\label{eq1.22A}
K_{\mu}(z)\quad=\quad\frac{1 + |\mu (z)|}{1 - |\mu\,(z)|}\,,
\end{equation}
$\mu(z)=\mu_f(z)=f_{\overline{z}}/f_z,$ при $f_z \ne 0$ и $\mu(z)=0$
в противном случае.}
\end{theorem}

\medskip
\begin{proof}
Заметим, что $f=\varphi\circ g,$ $g$ -- некоторый гомеоморфизм, а
$\varphi$ -- аналитическая функция (см. \cite[п.~5 (III),
гл.~V]{St}). Следовательно, отображение $f$ дифференцируемо почти
всюду (см., напр., \cite[Теорема~3.1,\,\S\,3, гл.~III]{LV}). Пусть
$B$ -- борелево множество всех точек $z\in D,$ где $f$ имеет полный
дифференциал $f^{\,\prime}(z)$ и $J(z, f)\ne 0$. Заметим, что $B$
может быть представлено в виде не более, чем счётного объединения
борелевских множеств $B_l$, $l=1,2,\ldots\,,$ таких что
$f_l=f|_{B_l}$ являются билипшецевыми гомеоморфизмами (см.
\cite[пункты~3.2.2, 3.1.4 и 3.1.8]{Fe}). Без ограничения общности,
мы можем считать, что множества $B_l$ попарно не пересекаются.
Обозначим также символом $B_*$ множество всех точек $z\in D$ где $f$
имеет полный дифференциал, однако, $f^{\,\prime}(z)=0.$

Поскольку $f$ -- конечного искажения, $f^{\,\prime}(z)=0$ для почти
всех точек $z,$ где $J(z, f)=0.$ Таким образом, согласно построению
и учитывая сказанное, множество $B_0:=D\setminus \left(B\bigcup
B_*\right)$ имеет нулевую меру Лебега. Следовательно, по
\cite[теорема~9.1]{MRSY}, ${\mathcal H}^{\,1}(B_0\cap S_r)=0$ для
почти всех окружностей $S_r:=S(z_0,r)$ с центром в точке
$z_0\in\overline{D},$ где ${\mathcal H}^{\,1},$ как обычно, линейная
мера Хаусдорфа, а $"$почти всех$"$ следует понимать в смысле модуля
семейств кривых. По лемме \ref{lem8.2.11} также и ${\mathcal
H}^{\,1}(B_0\cap S_r)=0$ для почти всех $r\in {\Bbb R},$

\medskip
Рассмотрим разбиение множества $D_*:=B(z_0, \varepsilon_0)\cap
D\setminus\{z_0\},$ $0<\varepsilon_0<d_0=\sup\limits_{z\in
D}|z-z_0|,$ на счётное число кольцевых сегментов $A_k,$
$k=1,2,\ldots\,.$ Пусть $\varphi_k$ -- вспомогательная
квазиизометрия, отображающая $A_k$ на прямоугольник
$\widetilde{A_k}$ такой, что дуги окружностей отображаются на
отрезки прямых. (Например, можно взять в качестве
$\varphi_k(\omega)=\log(\omega-z_0),$ $\omega\in A_k$). Рассмотрим
семейство отображений $g_k=f\circ\varphi_k^{\,-1},$
$g_k:\widetilde{A_k}\rightarrow {\Bbb C}.$ Заметим, что $g_k\in
W^{1, 1}_{\rm loc}$ (см. \cite[разд.~1.1.7]{Ma}), откуда, в
частности $g_k\in ACL$ (см. \cite[теоремы~1 и 2, п.~1.1.3,
$\S\,$1.1, гл.~I]{Ma}). Поскольку абсолютная непрерывность на
фиксированном отрезке влечёт $N$-свойство относительно линейной меры
Лебега (см. \cite[разд.~2.10.13]{Fe}), мы будем иметь, что
${\mathcal H}^{\,1}((g_k\circ\varphi_k)(B_0\cap A_k\cap
S_r))={\mathcal H}^{\,1}(f(B_0\cap A_k\cap S_r))=0$ и, значит, ввиду
полуаддитивности меры Хаусдорфа, ${\mathcal H}^{\,1}(f(B_0\cap
S_r))=0$ для почти всех $r\in {\Bbb R}.$

\medskip
Далее, покажем что ${\mathcal H}^{\,1}(f(B_*\cap S_r))=0$ для почти
всех $r\in {\Bbb R}.$ Действительно, пусть $\varphi_k,$ $g_k$ и
$A_k$ такие, как определено выше, $A_k=\{z\in {\Bbb C}:
z-z_0=re^{i\varphi}, r\in (r_{k-1}, r_{k}), \varphi\in (\psi_{k-1},
\psi_k)\},$ и пусть $S_k(r)$ часть сферы $S(z_0, r),$ принадлежащая
сферическому сегменту $A_k,$ т.е., $S_k(r)=\{z\in {\Bbb C}:
z-z_0=re^{i\varphi}, \varphi\in (\psi_{k-1}, \psi_k)\}.$ По
построению, $\varphi_k$ отображает $S_k(r)$ на сегмент $I(k,
r)=\{z\in {\Bbb C}: z=\log r+it, t\in (\psi_{k-1}, \psi_k).$
Применяя \cite[Теорема~3.2.5]{Fe}, мы получим, что
$${\mathcal H}^{\,1}(g_k(\varphi_k(B_*\cap S_k(r))))\le
\int\limits_{g_k(\varphi_k(B_*\cap S_k(r)))}N(y, g_k,
\varphi_k(B_*\cap S_k(r)))d{\mathcal H}^{1}y=$$
$$=\int\limits_{\varphi_k(B_*\cap
S_k(r))}|g_k^{\,\prime}(r+te)|dt=0$$
для почти всех $r\in (r_{k-1}, r_{k}).$ Из сказанного выше следует,
что ${\mathcal H}^{\,1}(f(B_*\cap S_k(r)))=0$ для почти всех $r\in
(r_{k-1}, r_{k}).$ Ввиду полуаддитивности хаусдорфовой меры,
${\mathcal H}^{\,1}(f(B_*\cap S_r))=0$ для почти всех $r\in {\Bbb
R},$ что и требовалось установить.

\medskip
Пусть $\Gamma$ -- семейство всех пересечений окружностей $S_r$,
$r\in(\varepsilon,\varepsilon_0)$,
$\varepsilon_0<d_0=\sup\limits_{z\in D}\,|z-z_0|,$ с областью $D.$
Для заданной допустимой функции $\rho_*\in{\rm adm}\,f(\Gamma),$
$\rho_*\equiv 0$ вне $f(D)$, полагаем $\rho\equiv 0$ вне $D$ и на
$B_0,$ и
$$\rho(z)\colon=\rho_*(f(z))\Vert f^{\,\prime}(z)\Vert \qquad{\rm for}\ z\in D\setminus B_0\,.$$
Для фиксированного множества $D_{r}^{\,*}\in f(\Gamma),$
$D_{r}^{\,*}=f(S_r\cap D),$ заметим, что
$$D_{r}^{\,*}=\bigcup\limits_{i=0}^{\infty} f(S_r\cap
B_i)\bigcup f(S_r\cap B_*)\,,$$
и, следовательно, для почти всех $r\in (0, \varepsilon_0)$
\begin{equation}\label{eq10a}
1\le \int\limits_{D^{\,*}_r}\rho_*(y)d{\mathcal A_*}=
\sum\limits_{i=0}^{\infty} \int \limits_{f(S_r\cap B_i)} N (y,
S_r\cap B_i)\rho_*(y) d{\mathcal H}^{1}y +
\end{equation}
$$+\int\limits_{f(S_r\cap B_*)}N (y,
S_r\cap B_*) \rho_*(y) d{\mathcal H}^{\,1}y\,.$$ Учитывая доказанное
выше, из (\ref{eq10a}) мы получаем, что
\begin{equation}\label{eq11A}
1\le \int\limits_{D^{\,*}_r}\rho_*(y)d{\mathcal A_*}=
\sum\limits_{i=1}^{\infty} \int \limits_{f(S_r\cap B_i)}N (y,
S_r\cap B_i) \rho_*(y) d{\mathcal H}^{\,1}y
\end{equation}
для почти всех $r\in (0, \varepsilon_0).$
Рассуждая покусочно на $B_i$, $i=1,2,\ldots$, согласно
\cite[пункт~1.7.6 и теорема~3.2.5]{Fe} мы получаем, что
$$\int\limits_{B_i\cap S_r}\rho\,d{\mathcal A}=
\int\limits_{B_i\cap S_r}\rho_*(f(z))\Vert
f^{\,\prime}(z)\Vert\,d{\mathcal A}=$$
$$=\int\limits_{B_i\cap S_r}\rho_*(f(z))\cdot\frac{\Vert
f^{\,\prime}(z)\Vert}{\frac{d{\mathcal A_*}}{d{\mathcal A}}}\cdot
\frac{d{\mathcal A_*}}{d{\mathcal A}}\,d{\mathcal A}\ge
\int\limits_{B_i\cap S_r}\rho_*(f(z))\cdot \frac{d{\mathcal
A_*}}{d{\mathcal A}}\,d{\mathcal A}=$$
\begin{equation}\label{eq12A}
=\int\limits_{f(B_i\cap S_r)}\rho_{*}\,d{\mathcal A_*}
\end{equation} для почти всех $r\in (0, \varepsilon_0).$
Из (\ref{eq11A}) и (\ref{eq12A}) с учётом леммы \ref{lem8.2.11}
следует, что $\rho\in{\rm{ext\,adm}}\,\Gamma.$

\medskip
Используя замену переменных на каждом $B_l$, $l=1,2,\ldots$ (см.,
напр., \cite[теорема~3.2.5]{Fe}), а также свойство счётной
аддитивности интеграла Лебега, получаем оценку
$$\int\limits_{D}\frac{\rho(z)}{K_{\mu}(z)}\,dm(z)\le
\int\limits_{f(D)}N(f, D)\rho_*(y)\, dm(y)\,,$$
что и завершает доказательство.
\end{proof} $\Box$

\section{Основные результаты}
\medskip{} Обозначим символом ${\frak F_{Q, E}}(x_0)$ семейство
всех открытых дискретных кольцевых $Q$-отоб\-ра\-же\-ний
$f:D\,\rightarrow\,\overline{{\Bbb R}^n}\setminus E$ в точке $x_0\in
D,$ а через ${\frak F_{Q, E}}(D)$ -- семейство открытых дискретных
кольцевых $Q$-отоб\-ра\-же\-ний $f:D\rightarrow\,\overline{{\Bbb
R}^n}\setminus E,$ при надлежащих классу ${\frak F_{Q, E}}(x_0)$ в
каждой точке $x_0\in D.$ Следующий результат был доказан в работе
\cite{Sev$_1$} (см. также статью \cite{RS}, где изложен аналогичный
случай гомеоморфизмов).

\medskip
\begin{proposition}\label{theor4*!} {\sl\, Пусть $E\subset\overline{{\Bbb R}^n}$ --
компактное множество положительной ёмкости. Тогда:

I. Семейство отображений ${\frak F_{Q, E}}(x_0)$ равностепенно
непрерывно в точке $x_0\in D,$ как только выполнено одно из
следующих условий:

1) $Q\in FMO(x_0);$

2) $q_{x_0}(r)=\,O\left({\left[
\log{\frac{1}{r}}\right]}^{n-1}\right)$ при $r\rightarrow 0,$ где
$q_{x_0}(r)$ определено соотношением (\ref{eq9b});

2) при некотором $\varepsilon_0=\varepsilon(x_0)>0,$
$\varepsilon(x_0)<{\rm dist\,}(x_0,\partial D),$ выполнено условие
(\ref{eq9}).

II. Семейство отображений ${\frak F_{Q, E}}(D)$ нормально, как
только хотя бы одно из условий 1)--3) выполнено в каждой точке
области $D.$} \end{proposition}

\medskip{} {\bf Доказательство} основных результатов работы --
теоремы \ref{th1}, следствия \ref{cor2}, теоремы \ref{th2} и
следствия \ref{cor3} -- напрямую вытекает из леммы \ref{lem1},
теоремы \ref{thOS4.1} и предложения \ref{theor4*!}. $\Box$

\medskip{} Ниже будет показано, что при $n=2$ утверждения теоремы \ref{th1}
верны при значительно более слабых ограничениях. Следующий результат
является важным непосредственно и с точки зрения дальнейших
приложений.

\medskip
\begin{theorem}\label{th6}
{\sl\, Пусть $D\subset {\Bbb C},$ $z_0\in D$ и $Q:D\rightarrow {\Bbb
C}$ -- локально интегрируемая функция. Тогда каждое открытое
дискретное отображение $f:D\rightarrow {\Bbb C}$ конечного искажения
класса $W^{1, 1}_{\rm loc}$ является кольцевым $Q$-отображением в
точке $z_0$ при $Q(z)=K_{\mu}(z),$ где $K_{\mu}(z)$ определено
соотношением (\ref{eq1.22A}).}
\end{theorem}

\medskip
{\it Доказательство} основано на теореме Стоилова о факторизации
(см. \cite[п.~5 (III), гл.~V]{St}. Согласно указанной теореме,
$f\,=\,\varphi\circ g,$
где $g$ -- некоторый гомеоморфизм, а $\varphi$ -- аналитическая
функция. Отметим, что множество точек ветвления $B_{\varphi}\subset
g(D)$ функции $\varphi$ состоит только из изолированных точек (см.
\cite[пункты 5 и 6 (II), гл.~V]{St}). Следовательно,
$g(z)=\varphi^{-1}\circ f$ локально, вне множества
$g^{-1}\left(B_{\varphi}\right).$ Ясно, что множество
$g^{-1}\left(B_{\varphi}\right)$ также состоит из изолированных
точек, следовательно, $g\in ACL(D)$ как композиция аналитической
функции $\varphi^{-1}$ и отображения $f\in W_{loc}^{1,1}(D).$

Покажем, что $g\in W_{loc}^{1,1}(D).$ Пусть далее $\mu_f(z)$
означает комплексную дилатацию функции $f(z),$ а $\mu_g(z)$ --
комплексную дилатацию $g.$ Согласно \cite[(1), п.~C, гл.~I]{A} для
почти всех $z\in D$ получаем:
$$f_z=\varphi_z(g(z))g_z,\qquad f_{\overline{z}}=\varphi_z(g(z))g_{\overline{z}},$$
$$\mu_f(z)=\mu_g(z)=:\mu(z), \quad K_{\mu_f}(z)=K_{\mu_g}(z):=K_{\mu}(z)=\frac{1+|\mu|}{1-|\mu|}\,.$$
Таким образом,  $K_{\mu}(z)\in L_{loc}^1(D).$ Поскольку $f$ --
конечного искажения, $g$ также конечного искажения и при почти всех
$z\in D$ выполнены соотношения
$$|\partial g|\le |\partial g|+ |\overline{\partial} g|= K^{1/2}_{\mu}(z)J^{1/2}(f, z)\,,$$
откуда по неравенству Гёльдера $|\partial g|\in L_{loc}^1 (D)$ и
$|\overline{\partial} g|\in L_{loc}^1 (D).$ Следовательно, $g\in
W_{loc}^{1,1}(D)$ и $g$ имеет конечное искажение. По теореме
\ref{thOS4.2} отображение $g$ является нижним $Q$-отображением в
точке $z_0$ при $Q=K_{\mu}(z)$ и, следовательно, $g$ -- кольцевое
$Q$-отображение в точке $z_0$ ввиду леммы \ref{lem1}. Завершает
доказательство применение неравенства Полецкого (см.
\cite[теорема~1]{Pol}), согласно которому аналитическая функция
$\varphi$ удовлетворяет условию $M(\varphi(\Gamma))\le M(\Gamma)$
для любого семейства кривых $\Gamma$ в $g(D)$ и, таким образом,
ввиду представления $f=\varphi\circ g$ отображение $f$ также
является кольцевым $Q$-отображением в точке $z_0.$ $\Box$

\medskip{} Ещё несколько результатов, связанных с равностепенной
непрерывностью классов Соболева на плоскости, могут быть
непосредственно получены из теоремы \ref{th6} и предложения
\ref{theor4*!}. Для заданных компактного множества $E\subset{\Bbb
C},$ области $D\subset {\Bbb C}$ и измеримой по Лебегу функции
$Q:D\rightarrow [1, \infty]$ обозначим символом $\frak{B}_{Q, E}$
семейство всех открытых дискретных отображений $f:D\rightarrow {\Bbb
C}\setminus E$ класса $W^{1, 1}_{loc},$ имеющих конечное искажение,
таких что $K_{\mu}(z)\le Q(z)$ почти всюду (где $K_{\mu}(z)$
определено соотношением (\ref{eq1.22A})). Справедливы следующие
результаты.

\medskip
\begin{theorem}\label{th3}
{\sl\, Семейство отображений $\frak{B}_{Q, E}$ является
равностепенно непрерывным в некоторой фиксированной точке $x_0\in
D,$ если ${\rm cap\,}E>0,$ $Q\in L_{loc}^1(D)$ и
при некотором $\varepsilon_0>0,$ $\varepsilon_0<{\rm dist}(x_0,
\partial D),$ выполнено условие расходимости интеграла (\ref{eq9})
(где $n=2$);
здесь, как обычно, $q_{x_0}(r)$ задаётся соотношением (\ref{eq9b}).
В частности, заключение теоремы \ref{th3} является верным, если
$q_{x_0}(r)=\,O\left({\left[ \log{\frac{1}{r}}\right]}\right)$ при
$r\rightarrow 0.$}
\end{theorem}

\medskip
Из критерия Арцела-Асколи вытекает следующее

\medskip
\begin{corollary}\label{cor4}{\sl\,
В условиях теоремы \ref{th3} семейство отображений $\frak{B}_{Q, E}$
является нормальным семейством отображений, как только условие
(\ref{eq9}) выполнено в каждой точке $x_0$ области $D.$}
\end{corollary}

\medskip
\begin{theorem}\label{th4}
{\sl\, Семейство отображений $\frak{B}_{Q, E}$ является
равностепенно непрерывным в некоторой фиксированной точке $x_0\in
D,$ если выполнено условие: $Q\in FMO(x_0).$}
\end{theorem}

\medskip
На основании критерия Арцела--Асколи вытекает следующее

\medskip
\begin{corollary}\label{cor5}{\sl\,
В условиях теоремы \ref{th4}, семейство отображений $\frak{B}_{Q,
E}$ является нормальным семейством отображений, как только условие
$Q\in FMO(x_0)$ выполнено в каждой точке $x_0$ области $D.$}
\end{corollary}

\section{Некоторые примеры}

\medskip{} Прежде всего, приведём ниже некоторое правило подсчёта
для величины $K_O(x, f),$ определённой в соотношении
(\ref{eq0.1.1A}). Предположим, что отображение $f:D\rightarrow {\Bbb
R}^n$ дифференцируемо в точке $x_0\in D$ и матрица Якоби
$f^{\,\prime}(x_0)$ невырождена, $J(x_0, f)={\rm
det\,}f^{\,\prime}(x_0)\ne 0.$ Тогда найдутся системы векторов
$e_1,\ldots, e_n$ и $\widetilde{e_1},\ldots,\widetilde{e_n}$ и
положительные числа
$$\lambda_1(x_0),\ldots,\lambda_n(x_0),\,\,\lambda_1(x_0)\le\ldots\le\lambda_n(x_0)\,,$$
такие что $f^{\,\prime}(x_0)e_i=\lambda_i(x_0)\widetilde{e_i}$ (см.
\cite[теорема~2.1 гл. I]{Re}), при этом,
$$|J(x_0, f)|=\lambda_1(x_0)\ldots\lambda_n(x_0),\quad \Vert f^{\,\prime}(x_0)\Vert
=\lambda_n(x_0)\,, $$
%
$$K_O(x_0,f)=\frac{\lambda^n_n(x_0)}{\lambda_1(x_0)\ldots\lambda_n(x_0)}\,,$$
%
см. \cite[соотношение~(2.5), разд.~2.1, гл.~I]{Re}.  Числа
$\lambda_1(x_0),\ldots\lambda_n(x_0)$ называются {\it главными
значениями}, а вектора $e_1,\ldots, e_n$ и
$\widetilde{e_1},\ldots,\widetilde{e_n}$ -- {\it главными векторами
} отображения $f^{\,\prime}(x_0).$

\medskip Заметим, что производная
$\frac{\partial f}{\partial e}(x_0)=\lim\limits_{t\rightarrow
+0}\frac{f(x_0+te)-f(x_0)}{t}$ отображения $f$ по направлению $e\in
{\Bbb S}^{n-1}$ в точке его дифференцируемости $x_0$ может быть
вычислена по правилу: $\frac{\partial f}{\partial
e}(x_0)=f^{\,\prime}(x_0)e.$ Таким образом, путём прямого вычисления
можно убедиться в справедливости следующего утверждения.

\medskip
\begin{proposition}\label{pr1C}
{\sl\, Пусть отображение $f: B(0, p)\rightarrow {\Bbb R}^n$ имеет
вид
\begin{equation}\label{eq8C}
f(x)=\frac{x}{|x|}\rho(|x|)\,,
\end{equation}
где функция $\rho(t):(0, p)\rightarrow {\Bbb R}$ непрерывна и
дифференцируема почти всюду. Тогда $f$ также дифференцируемо почти
всюду, при этом, в каждой точке $x_0$ дифференцируемости отображения
$f$ в качестве главных векторов $e_{i_1},\ldots, e_{i_n}$ и
$\widetilde{e_{i_1}},\ldots, \widetilde{e_{i_n}}$ можно взять
$(n-1)$ линейно независимых касательных векторов к сфере $S(0, r)$ в
точке $x_0,$ где $|x_0|=r,$ и один ортогональный к ним вектор в
указанной точке.

Соответствующие главные растяжения (называемые, соответственно, {\it
касательными растяжениями} и {\it радиальным растяжением}) равны
$\lambda_{\tau}(x_0):=\lambda_{i_1}(x_0)=\ldots=\lambda_{i_{n-1}}(x_0)=\frac{\rho(r)}{r}$
и $\lambda_{r}(x_0):=\lambda_{i_n}=\rho^{\,\prime}(r),$
соответственно.}
\end{proposition}

\medskip
Отметим, что для главных растяжений $\lambda_{i_k},$ $k\in 1,2,
\ldots, n,$ мы намеренно использовали двойную индексацию, поскольку,
как мы условились выше, конечную последовательность $\lambda_i,$
$i\in 1,2,\ldots, n$ мы предполагаем возрастающей по $i:$
$\lambda_1\le\lambda_2\le\ldots\le\lambda_n.$ Естественно, что в
фиксированной точке $x_0$ радиальные растяжения
$\lambda_{i_1}(x_0)=\ldots=\lambda_{i_{n-1}}(x_0)=\frac{\rho(r)}{r}$
могут быть не больше касательного растяжения
$\lambda_{i_n}=\rho^{\,\prime}(r),$ и наоборот.

\medskip{}
Отметим, что ограничения на функцию $Q,$ присутствующие в
формулировках основных результатов настоящей работы, нельзя, вообще
говоря, заменить условием $Q\in L^p$ ни для какого (сколь угодно
большого) $p>0$ и для любой неубывающей функции $\varphi(t).$ Для
простоты рассмотрим случай, когда $D={\Bbb B}^n,$ $n\ge 2.$ Имеет
место следующее утверждение.

\medskip
\begin{theorem}\label{th3.10.1}{\sl\,
Пусть $\varphi:[0,\infty)\rightarrow[0,\infty)$ -- произвольная
неубывающая функция. Для каждого $p\ge 1$ существуют функция
$Q:{\Bbb B}^n\rightarrow [1, \infty],$ $Q(x)\in L^p({\Bbb B}^n)$ и
равномерно ограниченная последовательность гомеоморфизмов $g_m:{\Bbb
B}^n\rightarrow {\Bbb R}^n,$ $g_m\in W_{loc}^{1, \varphi}({\Bbb
B}^n),$ имеющих конечное искажение, таких что $K_O^{n-1}(x, g_m)\le
Q(x),$ при этом, семейство $\left\{g_m(x)\right\}_{m=1}^{\infty}$ не
является равностепенно непрерывным в точке $x_0=0.$}
\end{theorem}

\medskip
\begin{proof} Рассмотрим следующий пример.
Зафиксируем числа $p\ge 1$ и $\alpha\in \left(0, n/p(n-1)\right).$
Можно считать, что $\alpha<1$ в силу произвольности выбора $p.$
Зададим последовательность гомеоморфизмов $g_m: {\Bbb
B}^n\rightarrow {\Bbb R}^n$ следующим образом:
$$ g_m(x)\,=\,\left
\{\begin{array}{rr} \frac{1+|x|^{\alpha}}{|x|}\cdot x\,, & 1/m\le|x|<1, \\
\frac{1+(1/m)^{\alpha}}{(1/m)}\cdot x\,, & 0<|x|< 1/m \ .
\end{array}\right.
$$
Заметим, что каждое отображение $g_m$ переводит шар $D={\Bbb B}^n$ в
шар $D^{\,\prime}=B(0,2)$ и что последовательность $g_m$ постоянна
при $|x|\ge 1/m,$ а именно, $g_m(x)\equiv g(x)$ при всех $x:\
\frac{1}{m}<|x|< 1,$ $m=1,2\ldots\,,$ где
$g(x)=\frac{1+|x|^{\alpha}}{|x|}\cdot x.$

Заметим, что $g_m\in ACL({\Bbb B}^n).$ Действительно, отображения
$g_m^{(1)}(x)=\frac{1+(1/m)^{\alpha}}{(1/m)}\cdot x,$
$m=1,2,\ldots,$ являются отображениями класса $C^1,$ скажем, в шаре
$B(0, 1/m+\varepsilon)$ при малых $\varepsilon>0,$ а отображения
$g_m^{(2)}(x)=\frac{1+|x|^{\alpha}}{|x|}\cdot x$ -- отображениями
класса $C^1,$ скажем, в кольце $$A(1/m-\varepsilon, 1,
0)=\left\{x\in {\Bbb R}^n: 1/m-\varepsilon<|x|<1\right\}$$ при малых
$\varepsilon>0.$ Отсюда вытекает, что гомеоморфизмы $g_m$ являются
липшицевыми в ${\Bbb B}^n$ и, значит, $g_m\in ACL({\Bbb B}^n)$ (см.,
напр., \cite[разд.~5, с.~12]{Va}).

Далее, в каждой регулярной точке $x\in D$ отображения
$g_m:D\rightarrow {\Bbb R}^n$ вычислим внешнюю дилатацию отображения
$g_m$ в точке $x$ (см. соотношение (\ref{eq0.1.1A})). Поскольку
каждое $g_m$ имеет вид (\ref{eq8C}), согласно предложению \ref{pr1C}
получаем, что, во-первых, $K_O(x, g_m)=1$ при $x\in B(0, 1/m),$
во-вторых, при $1/m\le|x|< 1$ имеем, что
$\lambda_{\tau}(x)=\frac{|x|^{\alpha}+1}{|x|},$
$\lambda_{r}(x)=\alpha|x|^{\alpha-1},$ $\Vert
g_m^{\,\prime}(x)\Vert=\frac{|x|^{\alpha}+1}{|x|},$ $|J(x, g_m)|=
\left(\frac{|x|^{\alpha}+1}{|x|}\right)^{n-1}\cdot
\alpha|x|^{\alpha-1}$ и
$$K_O(x, g_m)=\left
\{\begin{array}{rr} \frac{1+|x|^{\,\alpha}}{\alpha
|x|^{\,\alpha}}\,, & 1/m\le|x|\le 1, \\
1\,,\qquad & 0<|x|< 1/m\,.
\end{array}\right.
$$
Заметим, что при каждом фиксированном $m\in {\Bbb N}$ и некотором
$c_m>0$ имеет место неравенство $\Vert g_m^{\,\prime}(x)\Vert\le
c_m,$ кроме того, нетрудно видеть, что $|\nabla g_m(x)|\le
n^{1/2}\cdot\Vert g_m^{\,\prime}(x)\Vert$ при почти всех $x\in {\Bbb
B}^n.$ Тогда ввиду неубывания функции $\varphi$
$$\int\limits_{{\Bbb B}^n} \varphi(|\nabla g_m(x)|)dm(x)\le\varphi(n^{1/2}c_m)\cdot m({\Bbb B}^n)<\infty\,,$$
т.е., $g_m\in W^{1, \varphi}({\Bbb B}^n).$
Заметим, что отображения $g_m$ имеют конечное искажение, поскольку
их якобиан почти всюду не равен нулю; кроме того, $K_O^{n-1}(x,
f)\le Q(x),$ где $Q=\left(\frac{1+|x|^{\,\alpha}}{\alpha
|x|^{\,\alpha}}\right)^{n-1},$ и
$$Q(x)\le \frac{C}{|x|^{\alpha(n-1)}}\,,\quad C:=\left(\frac{2}{\alpha}\right)^{n-1}\,.$$
Таким образом, получаем: $$\int\limits_{{\Bbb
B}^n}\left(Q(x)\right)^p dm(x)\le C^p \int\limits_{{\Bbb
B}^n}\frac{dm(x)}{|x|^{p\alpha(n-1)}}=$$
\begin{equation}\label{eq2.3A}=C^p\int\limits_0^1\int\limits_{S(0,
r)}\frac{d{\mathcal{A}}}{|x|^{p\alpha(n-1)}}\,dr=\omega_{n-1}C^p
\int\limits_0^1\frac{dr}{r^{(n-1)(p\alpha-1)}}\,.\end{equation}
Хорошо известно, что интеграл
$I:=\int\limits_0^1\frac{dr}{r^{\beta}}$ сходится при $\beta<1.$
Таким образом, интеграл в правой части соотношения (\ref{eq2.3A})
сходится, поскольку показатель степени $\beta:=(n-1)(p\alpha-1)$
удовлетворяет условию $\beta<1$ при $\alpha\in (0, n/p(n-1)).$

Отсюда вытекает, что $Q(x)\in L^p({\Bbb B}^n).$ С другой стороны,
легко видеть, что
\begin{equation}\label{eq2!!!!!}
\lim\limits_{x\rightarrow 0} |g(x)|= 1\,,
\end{equation}
и $g$ отображает проколотый шар ${\Bbb B}^n\setminus\{ 0\}$ на
кольцо $1<|y|< 2.$ Тогда, ввиду (\ref{eq2!!!!!}), мы получаем, что
$$|g_m(x)|=|g(x)|\ge 1\qquad\qquad\forall\quad x:|x|\ge 1/m,\quad m=1,2,\ldots\,,$$
т.е., семейство $\{g_m\}_{m=1}^{\infty}$ не является равностепенно
непрерывным в нуле.
\end{proof}

\medskip{} Приведём ещё один интересный, на наш взгляд, пример,
касающийся выполнения условия (\ref{eq9}) в формулировках основных
утверждений работы. Хотя мы и не можем в буквальном смысле назвать
это условие необходимым и достаточным условием равностепенной
непрерывности соответствующего семейства отображений, условие
(\ref{eq9}), всё же, является условием $"$близким$"$ к необходимому
в следующем смысле.

\medskip
\begin{theorem}\label{th5} {\sl Пусть $\varphi:[0,\infty)\rightarrow[0,\infty)$ -- произвольная
неубывающая функция и $0<\varepsilon_0<1.$ Для каждой измеримой по
Лебегу функции $Q:{\Bbb B}^n\rightarrow [1, \infty],$ $Q\in
L_{loc}({\Bbb B}^n),$ такой, что
$\int\limits_{0}^{\varepsilon_0}\frac{dt}{tq_{0}^{\,\frac{1}{n-1}}(t)}<\infty,$
найдётся семейство равномерно ограниченных отображений $f_m\in
W_{loc}^{1, \varphi}({\Bbb B}^n)$ с конечным искажением со
следующими свойствами:

1) $K^{n-1}_O(x, f_m)\le \widetilde{Q}(x),$ где -- некоторая
измеримая по Лебегу функция, такая что
$\widetilde{q}_0(r):=\frac{1}{\omega_{n-1}r^{n-1}}\int\limits_{S(0,
r)}\widetilde{Q}(x)d{\mathcal H}^{n-1}=q_0(r)$ для почти всех $r\in
(0, 1);$

2) последовательность $f_m$ не является равностепенно непрерывной в
нуле.}
\end{theorem}

\begin{proof}
Определим последовательность отображений $f_m:{\Bbb B}^n\rightarrow
{\Bbb R}^n$ следующим образом:
$$f_m(x)=\frac{x}{|x|}\rho_m(|x|)\,,\qquad f_m(0):=0\,,$$
где
$$\rho_m(r)=
\exp\left\{-\int\limits_{r}^1\frac{dt}{tq_{0,
m}^{1/(n-1)}(t)}\right\}\,, \qquad q_{0,
m}(r):=\frac{1}{\omega_{n-1}r^{n-1}}\int\limits_{|x|=r}Q_m(x)\,d{\mathcal{A}}\,,
$$
$$Q_m(x)\quad=\quad \left \{\begin{array}{rr} Q(x) , & \
|x|> 1/m\ ,
\\ 1\ ,  &  |x|\le 1/m\,.
\end{array} \right.$$
Заметим, что $f_m\in ACL$ при любом $m\in {\Bbb N}$ и отображения
$f_m$ дифференцируемы почти всюду в ${\Bbb B}^n.$ Согласно
предложению \ref{pr1C}
$$\Vert
f_m^{\,\prime}(x)\Vert=\frac{\exp\left\{-\int\limits_{|x|}^1\frac{dt}{tq_{0,
m}^{1/(n-1)}(t)}\right\}}{|x|}\,, |J(x,
f_m)|=\frac{\exp\left\{-n\int \limits_{|x|}^1\frac{dt}{tq_{0,
m}^{1/(n-1)}(t)}\right\}}{|x|^nq_{0, m}^{1/(n-1)}(|x|)}\,.$$
Заметим, что $J(x, f_m)\ne 0$ при почти всех $x.$ Покажем теперь,
что $\varphi(|\nabla f_m(x)|)\in L^1({\Bbb B}^n).$ Используя теорему
Фубини, получаем:
$$\int\limits_{{\Bbb B}^n}\varphi(|\nabla f_m(x)|)\,dm(x)\le$$
$$\le\int\limits_{{\Bbb B}^n}\varphi(n^{1/2}\Vert f_m^{\,\prime}(x)\Vert)\,dm(x)=
\omega_{n-1}\int\limits_{0}^1
r^{n-1}\varphi\left(n^{1/2}\frac{\exp\left\{-\int\limits_{r}^1\frac{dt}{tq_{0,
m}^{1/(n-1)}(t)}\right\}}{r}\right)dr=$$
\begin{equation}\label{eq14A}
=\omega_{n-1}\left(\int\limits_{0}^{1/m}\psi(r)dr+\int\limits_{1/m}^1\psi(r)dr\right)
=\omega_{n-1}(I_1+I_2)\,,
\end{equation}
где $I_1:=\int\limits_{0}^{1/m}\psi(r)dr,$
$I_2:=\int\limits_{1/m}^1\psi(r)dr$ и
$\psi(r):=r^{n-1}\varphi\left(n^{1/2}\frac{\exp\left\{-\int\limits_{r}^1\frac{dt}{tq_{0,
m}^{1/(n-1)}(t)}\right\}}{r}\right).$ Заметим, прежде всего, что
$I_2\le \varphi(n^{1/2}m)\cdot \frac{m-1}{m}\le \varphi(n^{1/2}m).$
Далее, заметим, что $I_1\le
\int\limits_{0}^{1/m}r^{n-1}\varphi\left(n^{1/2}\frac{\exp\left\{-\int\limits_{r}^{1/m}\frac{dt}{tq_{0,
m}^{1/(n-1)}(t)}\right\}}{r}\right)dr\le\varphi(n^{1/2}m).$ В таком
случае, из соотношений (\ref{eq14A}) вытекает, что
$\int\limits_{{\Bbb B}^n}\varphi(|\nabla f_m(x)|)\,dm(x)\le
2\omega_{n-1}\cdot\varphi(n^{1/2}m),$ т.е., $\varphi(|\nabla
f_m(x)|)\in L^1({\Bbb B}^n).$

Заметим, что $K_O(x, f_m)=q_{0, m}^{1/(n-1)}(|x|)\le
q_{0}^{1/(n-1)}(|x|)$ и, значит, $K_O^{n-1}(x, f_m)\le q_{0}(|x|)$
при почти всех $x\in {\Bbb B}^n.$ Полагаем
$\widetilde{Q}(x):=q_{0}(|x|),$ тогда будем иметь, что
$\widetilde{q}_0(r)=q_0(r)$ для почти всех $r\in (0, 1).$

\medskip
Заметим, что $|f_m(x)|\le 1$ для всех $m\in {\Bbb N}$ и, таким
образом, семейство отображений $\{f_l(x)\}_{l=1}^{\infty}$
равномерно ограничено. Осталось показать, что построенная, таким
образом, последовательность отображений $f_m$ не является
равностепенно непрерывной в нуле. Для произвольной
последовательности $x_m$ такой, что $|x_m|=1/m,$ $m=1,2,\ldots,$
имеем $|f_m(x_m)|\ge \sigma,$ где $\sigma$ не зависит от $m.$
Окончательно, для некоторого числа $\sigma$ и произвольного элемента
последовательности $1/(m-1),$ $m=2,3,\ldots,$ найдётся $x_m\in {\Bbb
B}^n$ и элемент семейства отображений
$f_m\in\{f_l(x)\}_{l=1}^{\infty}$ такие, что $|x_m-0|<1/(m-1)$ и, в
то же время, $|f_m(x_m)-f_m(0)|\ge \sigma.$ Таким образом, семейство
отображений $\{f_l(x)\}_{l=1}^{\infty}$ не является равностепенно
непрерывным в нуле.
\end{proof}

\medskip
\noindent{{\bf Евгений Александрович Севостьянов} \\ Институт
прикладной математики и механики НАН Украины \\
83 114 Украина, г. Донецк, ул. Розы Люксембург, д. 74, \\
тел. +38 (066) 959 50 34 (моб.), +38 (062) 311 01 45 (раб.), e-mail:
brusin2006@rambler.ru, esevostyanov2009@mail.ru}
\end{document}